\documentclass[mathpazo]{cicp}

\usepackage{lipsum}
\usepackage{amsmath}
\usepackage{amsfonts}
\usepackage{graphicx}
\usepackage{epstopdf}
\usepackage{bm}
\usepackage{amsthm} 
\usepackage{amssymb}
\usepackage{algorithm}
\usepackage{algpseudocode}
\usepackage{listings}
\usepackage{booktabs}
\usepackage{epsfig}
\usepackage{color}
\usepackage{float}
\usepackage{footnote}
\usepackage{mathrsfs}
\usepackage{multirow}
\usepackage{siunitx}
\usepackage{subfig}
\usepackage{tikz}
\usetikzlibrary{shapes}
\usetikzlibrary{positioning}
\usetikzlibrary{circuits}
\usepackage{indentfirst}
\usepackage{hyperref}
\usetikzlibrary{positioning, shapes, calc, fit}
\ifpdf
\DeclareGraphicsExtensions{-eps-converted-to.pdf,.pdf,.png,.jpg}
\else
\DeclareGraphicsExtensions{-eps-converted-to.pdf}
\fi

\newcommand{\x}{\bm{\mathrm{x}}}
\newcommand{\bu}{\bm{\mathrm{u}}}

\usepackage{tikz}
\usetikzlibrary{positioning, shapes, arrows.meta, fit, calc}

\begin{document}
\title{Structured First-Layer Initialization Pre-Training Techniques to Accelerate Training Process Based on $\varepsilon$-Rank}


\author[Tang T et.~al.]{Tao Tang\affil{3},
      Jiang Yang\affil{1,2}\comma\corrauth, Yuxiang Zhao\affil{1} and Quanhui Zhu\affil{1}}
\address{\affilnum{1}\ Department of Mathematics, Southern University of Science and Technology, Shenzhen,China. \\
          \affilnum{2}\  SUSTech International Center for Mathematics, Shenzhen, China.\\
          \affilnum{3}\ School of Mathematics and Statistics, Guangzhou Nanfang College, Guangzhou, China.}
\emails{{\tt ttang@nfu.edu.cn} (T.~Tang), {\tt yangj7@sustech.edu.cn} (J.~Yang),
        {\tt 12131241@mail.sustech.edu.cn} (Y.~Zhao), {\tt 12131244@mail.sustech.edu.cn} (Q.~Zhu)}

\begin{abstract}
Training deep neural networks for scientific computing remains computationally expensive due to the slow formation of diverse feature representations in early training stages. Recent studies identify a staircase phenomenon in training dynamics, where loss decreases are closely correlated with increases in $\varepsilon$-rank, reflecting the effective number of linearly independent neuron functions. 
 Motivated by this observation, this work proposes a structured first-layer initialization (SFLI) pre-training method to enhance the diversity of neural features at initialization by constructing $\varepsilon$-linearly independent neurons in the input layer. We present systematic initialization schemes compatible with various activation functions and integrate the strategy into multiple neural architectures, including modified multi-layer perceptrons and physics-informed residual adaptive networks. Extensive numerical experiments on function approximation and PDE benchmarks, demonstrate that SFLI significantly improves the initial $\varepsilon$-rank, accelerates convergence, mitigates spectral bias, and enhances prediction accuracy. With the help of SILP, we only need to add one line of code to conventional existing algorithms.
\end{abstract}


\maketitle


\section{Introduction}
Neural networks have become a cornerstone of modern machine learning, achieving remarkable accuracy in both classical machine learning tasks and emerging areas such as scientific computing and partial differential equation (PDE) modeling. Despite their success, training such models remains computationally intensive, often requiring large-scale resources and prolonged optimization. This has motivated a broad range of efforts to accelerate the training process, including architectural innovations, advanced optimization techniques, and improvements in loss function design.

Among these efforts, particular attention has been paid to the learning dynamics of neural networks \cite{haninHowStartTraining2018, liCyclicalLearningRate2020,geigerLandscapeTrainingRegimes2021,nakhodnovLossFunctionDynamics2022,liVisualizingLossLandscape2018, zhiweiLossJumpLoss2024}.
A key observation in training dynamics is the so-called frequency principle (F-Principle) or spectral bias\cite{xuFrequencyPrincipleFourier2020,rahamanSpectralBiasNeural2019,xuOverviewFrequencyPrinciple2024,xuTrainingBehaviorDeep2019,zhangWhyShallowNetworks2023}, which states that neural networks tend to fit low-frequency components of the target function earlier in training. 
To address this, multiscale network designs like MscaleDNN \cite{liuMultiScaleDeepNeural2020,liSolvingClassMultiscale2024} have been proposed, which incorporate hierarchical frequency encodings. The performance of such architectures strongly depends on the behavior of the first hidden layer, particularly the diversity of its initial activations—a critical aspect that remains insufficiently studied.


A standard multi-layer perceptron (MLP)  can be formulated as follows:
\begin{equation}
	\left\{
	\begin{aligned}
		y_0 &= x, \\
		y_l &= H(y_{l-1};\theta_l),& l = 1,\cdots,L,\\
		y &= \beta\cdot y_L,
	\end{aligned}\right.
	\label{st::NN}
\end{equation}
where $x\in\mathbb{R}^d$ is the input, $y_l\in \mathbb{R}^{n_l}$ is the neurons of the $l$-th hidden layer, and $\beta\in\mathbb{R}^{n_L}$ is the coefficients in the output layer. Each layer mapping $H$ represents the structure of the hidden layer, and a fully connected layer can be explicitly given by  
$H(y_{l-1};\theta_l) = \sigma(W_ly_{l-1}+b_l)$, where $W_l\in \mathbb{R}^{n_{l+1}\times {n_l}}, b_l\in\mathbb{R}^{n_{l+1}} ,$ are trainable parameters. The activation function $\sigma$ is applied element-wise. The neurons $y_L(x)$ in the final hidden layer can be viewed as a collection of scalar neuron functions defined on the input domain, and the final output $y$ of the network is a linear combination of these functions. 
This perspective aligns with interpretations in the deep finite element method~\cite{xiongDeepFiniteElement2025} and the finite neuron method~\cite{xuFiniteNeuronMethod2020}, where $y_L(x)$ serves as a set of basis functions.

The training dynamics of deep neural networks have attracted significant interest, particularly in understanding how low-dimensional representations are formed during optimization \cite{papyanPrevalenceNeuralCollapse2020, achilleEmergenceInvarianceDisentanglement2018}. Recent work \cite{yangEffectiveRankStaircase2024} has identified \textit{staircase phenomenon}: \textbf{In training dynamics, the loss function often decreases rapidly along with a significant growth of linear independence of neuron functions.}
A consistently low $\varepsilon$-rank of $y_L$ during training indicates a lack of functional diversity among neurons, which in turn limits the expressive power of the network. Such limitations pose a critical challenge for problems in scientific computing, where resolving fine-scale or high-frequency structures is essential.

Motivated by these insights, this work introduces a novel pre-training strategy termed structured first-layer initialization (SFLI), designed to increase the $\varepsilon$-rank of neuron functions at initialization. By enforcing $\varepsilon$-linear independence through carefully constructed weights in the first hidden layer, the method enhances feature diversity and accelerates convergence, while preserving numerical stability and incurring no additional computational cost.

The main contributions of this work are summarized as follows:
\begin{itemize}
    \item[1.] We propose a novel structured first-layer initialization pre-training strategy, which enhances initial feature diversity of neural functions and accelerates the training process. The method is activation-function agnostic, compatible with a wide range of neural architectures, and incurs negligible computational overhead, making it broadly applicable and easy to integrate into existing frameworks.
    \item[2.] We demonstrate the effectiveness of structured first-layer initialization across various function approximation and PDE-solving benchmarks, where the method consistently improves prediction accuracy, convergence speed, and numerical stability under different training scenarios, while also mitigating spectral bias through improved feature diversity.
\end{itemize}

The remainder of the paper is organized as follows. Section 2 reviews the concept of $\varepsilon$-rank and introduces the staircase phenomenon. Section 3 presents the proposed structured first-layer initialization and discusses its implementation across various network architectures. Section 4 provides numerical results to illustrate the effectiveness of the method. Conclusions and future works are summarized in Section 5.

\section{Staircase Phenomenon}
\begin{figure}[htbp]
	\centering
	\includegraphics[width=\textwidth]{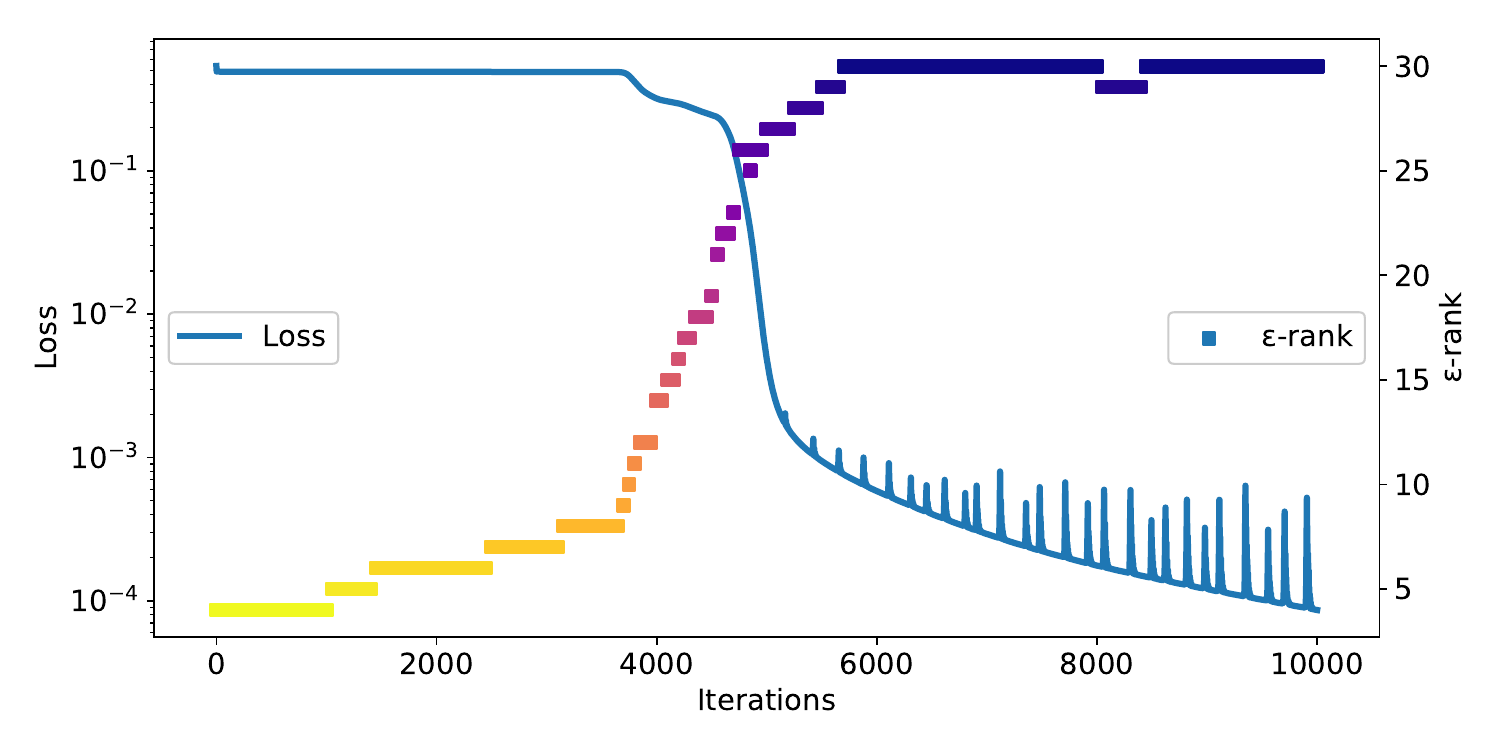}
	\caption{Staircase phenomenon in a regression task with target function $f(x)=\sin 20x$.}
	\label{fig::sp}
\end{figure}

In this section, we will first present a short recall on the $\varepsilon$-rank, which is introduced to quantitatively measure the linear independence of neuron functions, as formally defined below:
\begin{definition}[\cite{yangEffectiveRankStaircase2024}]
Let $u(x;\theta)$ denote a neural network with neuron functions $\{\phi_j\}_{j=1}^n$ in its final hidden layer. Define the Gram matrix 
 $\displaystyle (M_u)_{ij} = \int_\Omega \phi_i(x;\theta)\phi_j(x;\theta)\mathrm{d}x$. Then the $\varepsilon$-rank of $M_u$ is given by 
\begin{equation}
	r_\varepsilon(M_u) = \#\{\lambda \in \lambda(M_u) \mid \lambda > \varepsilon\},
\end{equation}
representing the number of eigenvalues $\lambda(M_u)$ exceeding a threshold $\varepsilon$. 
\end{definition}
Tracking $\varepsilon$-rank reveals the staircase phenomenon in training dynamics (Figure \ref{fig::sp}): the loss function $\mathcal{L}(u_n)$ decreases sharply when $r_\varepsilon(M_u)$ increases significantly. This behavior has been observed consistently across different tasks and architectures in \cite{yangEffectiveRankStaircase2024}.

\begin{figure}[htbp]
	\centering
	\includegraphics[width=0.9\textwidth]{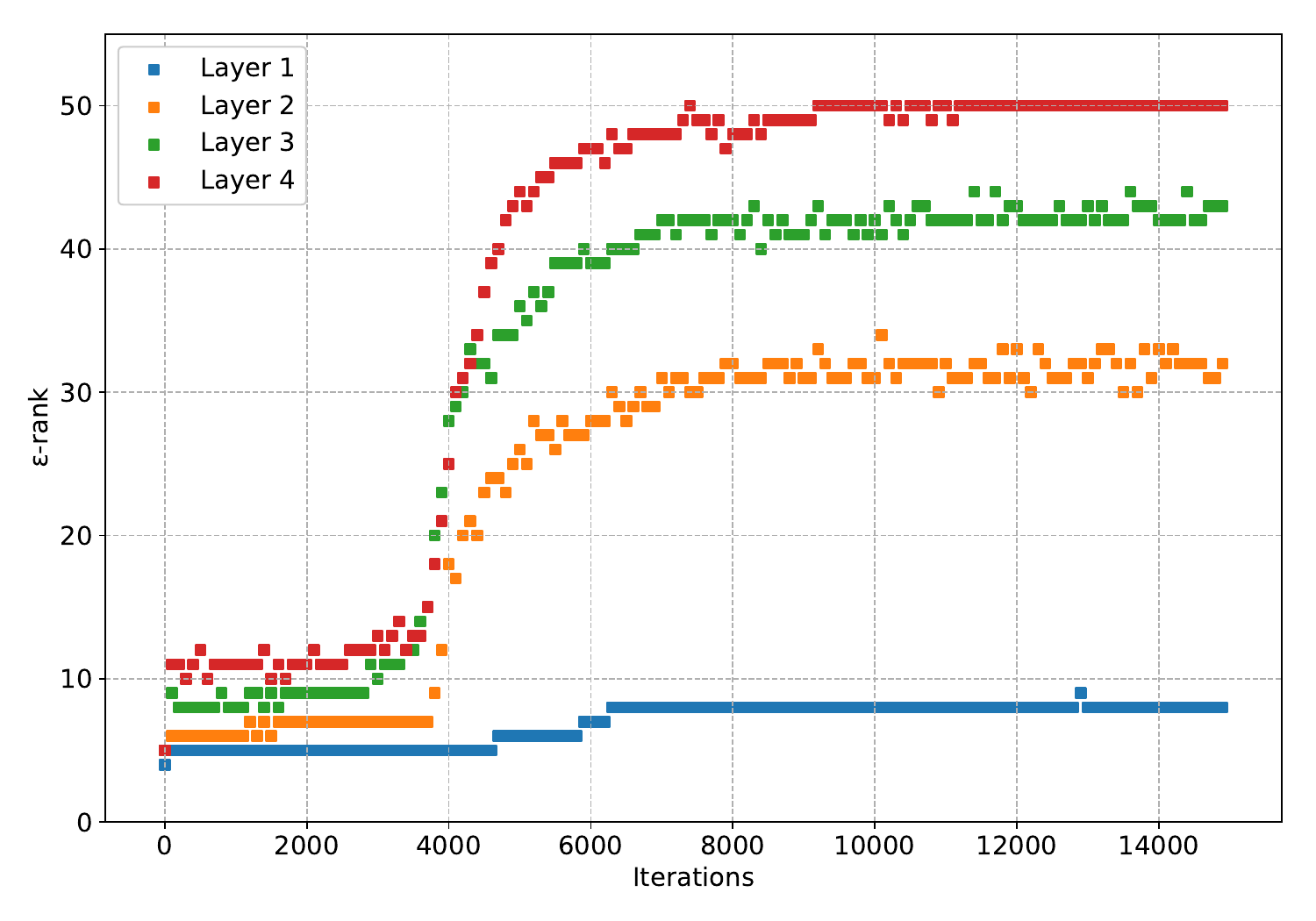}
	\caption{Layer-wise staircase phenomenon within the same neural network. The width is fixed $n_1=n_2=n_3=n_4=50$.}
	\label{fig::sp_depth}
\end{figure}

This phenomenon is theoretically justified by the following lower bound on the loss:
\begin{equation}
	\sqrt{\mathcal{L}(u_n)} \geq \frac{1}{C_s}\left(\sqrt{\text{dist}(u^*, \mathcal{F}_p)} - C(p+1)(n-p)^2\varepsilon\right),
\end{equation}
where $u_n\in\mathcal{F}_n$ is a neural network of width $n$, $p$ is the $\varepsilon$-rank of $u_n$, i.e., $r_\varepsilon(M_{u_n})=p$, and $\mathcal{F}_p$ denotes the function class with $p$ neurons in the last hidden layer. This inequality shows that meaningful reduction in loss requires an increase in the $\varepsilon$-rank.

Taken together, these observations suggest that a sufficiently large $\varepsilon$-rank is a prerequisite for significant loss decay.  Empirical evidence in \cite{yangEffectiveRankStaircase2024} further supports this claim: training curves often exhibit prolonged plateaus when $r_\varepsilon(M_u)$ remains stagnant, followed by decreases only after $\varepsilon$-rank jumps.

Standard initialization schemes, such as Xavier initialization \cite{glorotUnderstandingDifficultyTraining2010}, typically result in a low initial $\varepsilon$-rank, thereby requiring extended training to overcome the rank bottleneck. 
Figure \ref{fig::sp_depth} demonstrates the observation that within the same deep neural network, the $\varepsilon$-rank of the subsequent layer is higher than that of the previous layer. These insights have inspired pre-training strategies that directly construct $\varepsilon$-linearly independent neuron functions in the first hidden layer using deterministic weight initialization \cite{yangEffectiveRankStaircase2024}. This approach effectively elevates the initial $\varepsilon$-rank to $n$, bypassing the early-stage stagnation and accelerating convergence in practice.

This unified perspective, linking $\varepsilon$-rank to training dynamics, provides actionable guidelines for improving neural network performance. In the next section, we build on this foundation to propose a general, structured first-layer initialization strategy that systematically improves the initial representational capacity, with a focus on accelerating the training dynamics.

\section{Structured First-Layer Initialization}
\subsection{Pre-Activation Method}

Following the findings in \cite{yangEffectiveRankStaircase2024}, the first hidden layer plays a critical role in learning meaningful feature representations. Consider the neurons in the first hidden layer defined by:
\begin{equation}
	F(x) = \sigma(Wx+b), 
\end{equation}
where $x\in\Omega\subset \mathbb{R}^d$ is the input, $b=[b_1,\cdots, b_n]^\top\in\mathbb{R}^n$, $W=[w_1,\cdots,w_n]^\top\in\mathbb{R}^{n\times d}$ are trainable parameters, and $\sigma$ is an element-wise activation function. The vector-valued function $F(x) \in \mathbb{R}^n$ represents the output of the first hidden layer, where $F(x) := [F_1(x), \ldots, F_n(x)]^\top$ and each $F_i(x)$ is a scalar neuron function.

The objective of structured first-layer initialization is to construct a set of neuron functions that are approximately $\varepsilon$-linearly independent at initialization. 
This is achieved by strategically designing the weights and biases in the first hidden layer so that the neural functions $\{F_i(x)\}$ are well-separated in the input domain. The motivation stems from classical numerical methods such as finite element schemes, where basis functions are spatially localized and evenly distributed.

Specifically, we reparameterize each neuron function in the first hidden layer as:
\begin{equation}
    F_i(x;w_i,b_i) =\sigma(w_i\cdot x + b_i)= \sigma(\gamma(\alpha_i\cdot x + c_i)),\quad i=1,\cdots,n,
\end{equation}
where $w_i = \gamma \alpha_i,$ $b_i=\gamma c_i$. $c_i\ge0$ is the shift distributed uniformly in an interval with respect to $\Omega$, and $\alpha_i$ determines the orientation of the hyperplane. The parameter $\gamma>0$ controls the localization of the neural function, with a recommended practical choice:
\begin{equation}\label{choice}
    \gamma = C \cdot \frac{n^{1/d} - 1}{|\Omega|^{1/d}}, \quad C \in [0.5, 2].
\end{equation} 

\begin{figure}[htbp]
	\centering
	\includegraphics[width=\textwidth]{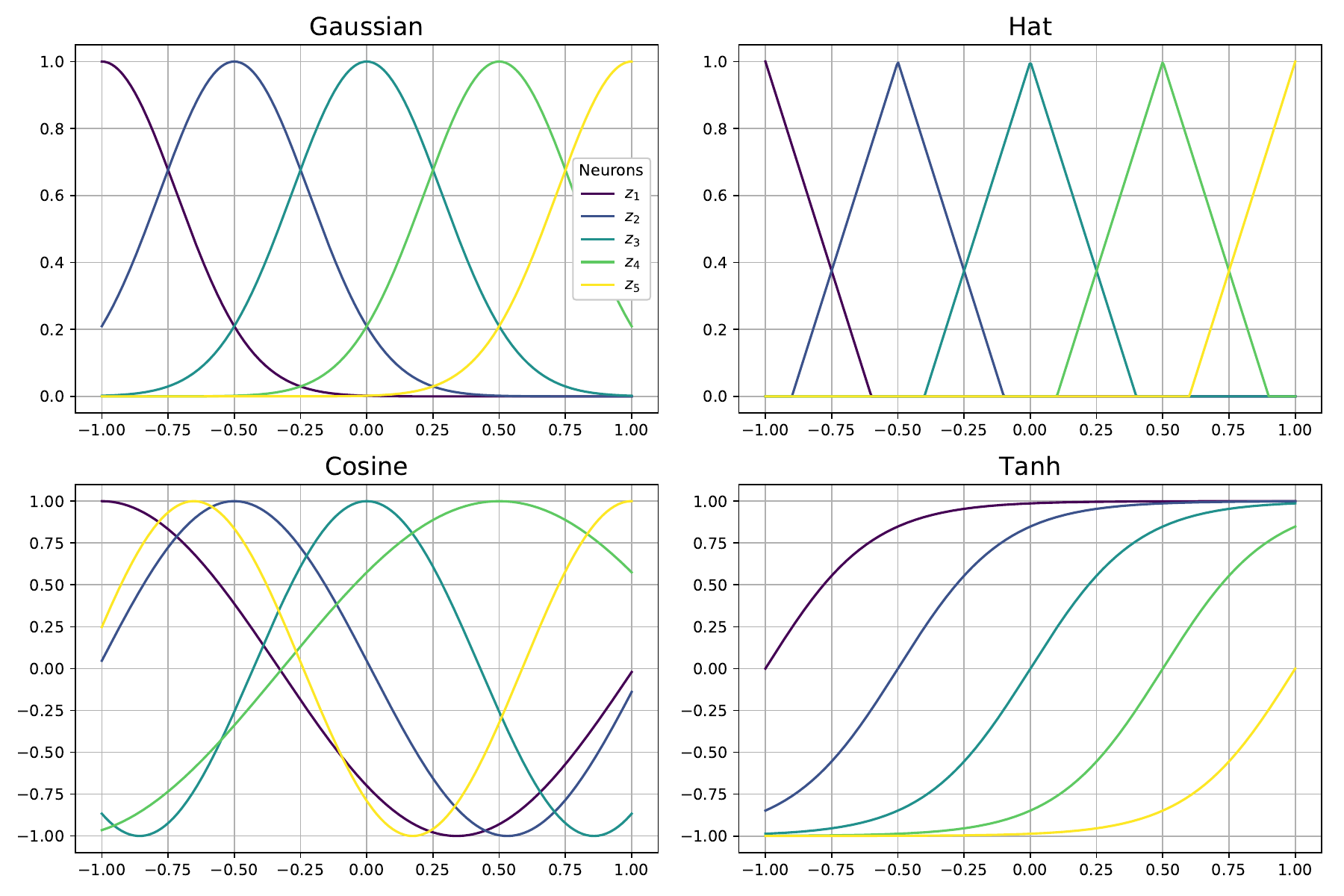}
	\caption{Visualization of the first-layer outputs using different activation functions under SFLI with $n=5$ neurons.}
	\label{fig::SFI_example}
\end{figure}

\begin{figure}[htbp]
\centering
\begin{tikzpicture}[
    every node/.style={font=\footnotesize},
    input/.style={circle, draw=blue!60, thick, minimum size=0.8cm},
    output/.style={circle, draw=blue!60, thick, minimum size=0.8cm},
    layer/.style={rectangle, draw=black!60, thick, rounded corners, minimum height=0.8cm, minimum width=2.6cm, fill=gray!10},
    neuron/.style={circle, draw=red!60, thick, minimum size=0.7cm},
    dashedbox/.style={draw=black!40, thick, rounded corners, dashed},
    connector/.style={-Stealth, thick}
  ]

\node[layer] (sfli-1) at (0, 1.5) {$F_1(x;w_1,b_1)$};
\node[layer] (sfli-2) [below=1 of sfli-1] {$F_2(x;w_2,b_2)$};
\node at ($(sfli-2) + (0,-1)$) (dots) {\vdots};
\node[layer] (sfli-4) [below=3.2cm of sfli-1] {$F_n(x;w_n,b_n)$};

\node[dashedbox, fit={(sfli-1)(sfli-2)(sfli-4)}] (sfliBox) {};

\node[input] (x) [left=1cm of sfliBox.west] {$x$};
\draw[connector] (x.east) -- (sfliBox.west);

\node[neuron] (hidden-11) [right=3.0cm of sfli-1] {};
\node[neuron] (hidden-12) [below=1cm of hidden-11] {};
\node at ($(hidden-12) + (0,-1.0cm)$) (hdots1) {\vdots};
\node[neuron] (hidden-13) [below=3.2cm of hidden-11] {};

\node[neuron] (hidden-21) [right=2.2cm of hidden-11] {};
\node[neuron] (hidden-22) [below=1cm of hidden-21] {};
\node at ($(hidden-22) + (0,-1.0cm)$) (hdots2) {\vdots};
\node[neuron] (hidden-23) [below=3.2cm of hidden-21] {};

\node at ($(hidden-11)!0.5!(hidden-21)$) {\(\cdots\)};

\node[output] (u) [right=2cm of $(hidden-21)!0.5!(hidden-23)$] {$u$};
\foreach \j in {21,22,23} {
  \draw[connector] (hidden-\j.east) -- (u.west);
}

\foreach \i in {1,2,4} {
  \foreach \j in {11,12,13} {
    \draw[connector] (sfli-\i.east) -- (hidden-\j.west);
  }
}

\node[dashedbox, fit={(hidden-11)(hidden-13)(hidden-21)(hidden-23)}] (hiddenBox) {};

\node[anchor=south, font=\small\bfseries] at (sfliBox.north) {SFLI Layer};
\node[anchor=south, font=\small\bfseries] at (hiddenBox.north) {Hidden Layers};

\end{tikzpicture}
\caption{Architecture of the neural network with Structured First-Layer Initialization (SFLI).}
\label{fig::sfli-architecture}
\end{figure}
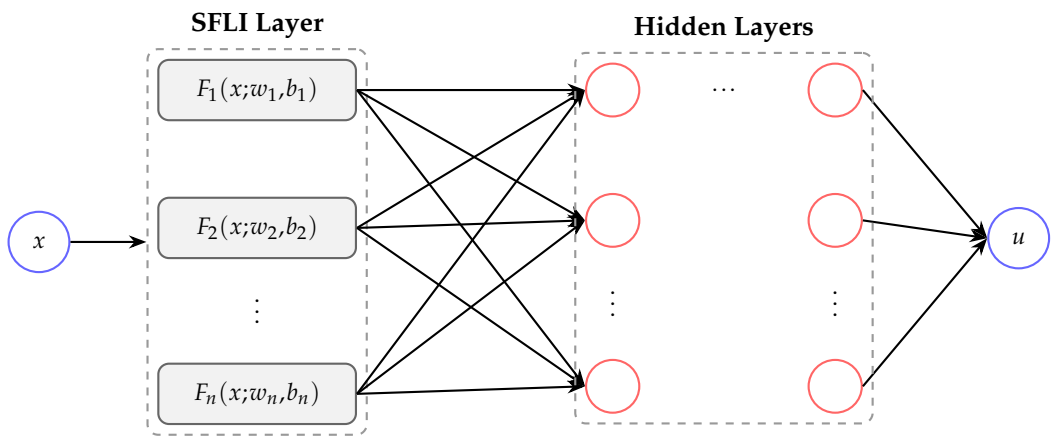

The overall structure of SFLI is illustrated in Figure~\ref{fig::sfli-architecture}. To provide an intuitive understanding, Figure~\ref{fig::SFI_example} visualizes the output behaviors of different activation functions under SFLI with $n=5$ neurons. Below, we list several typical activation functions along with their associated weight initialization and structural characteristics used in the SFLI method:
\begin{itemize}
	\item \textbf{Cosine:} $\sigma(x) = \cos(x)$. The weight vectors $\alpha_i$ are sampled independently from the standard normal distribution, i.e., $\alpha_i \sim  \mathcal{N}(\boldsymbol{0}, \boldsymbol{I}_d)
$.
	\item \textbf{Hyperbolic Tangent (Tanh):} $\sigma(x) = \tanh(x)$. The weights $\alpha_i$ are uniformly sampled from the unit sphere, i.e., $\alpha_i \in \mathbb{S}^d, \|\alpha_i\|_2 = 1$, to ensure directional diversity in the neuron responses.
	\item \textbf{Hat:} $\sigma(x) = \max(0, 1 - |x|)$. Similar to Tanh, the weight vectors are chosen from the unit sphere: $\alpha_i \in \mathbb{S}^d$. This ensures spatial coverage of the support regions of hat functions over the domain.
\item \textbf{Gaussian:} $\sigma(x) = e^{-|x|^2}$. 
 Unlike the previous activation functions, which are based on linear transformations $w_i \cdot x + b_i$, the Gaussian activation adopts a radial structure centered at specific locations in the input domain. Specifically, each neuron in the first hidden layer takes the form:
\begin{equation*}
	F_i(x;\gamma, c_i) = \exp\left(-\gamma^2 \|x - c_i\|^2\right),\quad i=1,\cdots,n,
\end{equation*}
where $\gamma > 0$ is a shape parameter that controls the sharpness of each Gaussian bump, and $c_i \in \mathbb{R}^d$ denotes the center of the $i$-th neuron uniformly distributed in $\Omega$. 

\end{itemize}

Importantly, SFLI modifies only the initialization of the first layer and does not change optimization, i.e., optimize $w_i, b_i$ directly rather than the decomposed parameters $\alpha_i, c_i$ and $\gamma$. The activation function and parameter settings are applied exclusively at initialization and do not affect subsequent layers, ensuring compatibility with standard training pipelines and incurring negligible computational cost.

The goals of the SFLI pre-training strategy can be summarized as follows:
\begin{itemize}
    \item Promote the diversity of the neurons over the input domain $\Omega$.
    \item Ensure that the $\varepsilon$-rank of the initial feature set  $\{F_i\}$ approaches $n$ with high probability.
    \item Maintain full compatibility with standard architectures without introducing additional training cost.
\end{itemize}

These principles underpin the design of the SFLI strategy and are further validated through empirical studies in the next section.

\subsection{Integrate with various neural network architecture}
The proposed structured first-layer initialization strategy can be seamlessly integrated into a wide range of neural network architectures, including Residual Neural Network (ResNet)\cite{heDeepResidualLearning2016}, modified MLP\cite{wangUnderstandingMitigatingGradient2021}, and Physics-Informed Residual AdapTivE Network (PirateNet)\cite{JMLR:v25:24-0313}. 
The modified MLP introduces a simple yet effective enhancement to the standard architecture, consistently outperforming the vanilla MLP in minimizing PDE residuals. PirateNet integrates adaptive residual connections to alleviate issues related to pathological initialization, enabling more stable and efficient training of physics-informed neural networks, especially in deeper architectures.
These models are representative of architectures commonly used in scientific computing and physics-informed learning.

As representative cases, we describe the network structures used in our numerical experiments. First, we introduce the modified MLP architecture:
\begin{equation}
	U = \sigma(W^{(u)}x+b^{(u)}),\quad V = \sigma(W^{(v)}x+b^{(v)}),
	\label{eq::encoder}
\end{equation}
followed by $L$ hidden layers:
\begin{equation}
	\begin{aligned}
		y_0 & = x,\\
		z_l &= \sigma(W_ly_{l-1}+b_l), \\
		y_l &= z_l\odot U + (1-z_l)\odot V,\quad l = 1,\dots,L,
	\end{aligned}
	\label{eq::modifiedMLP}
\end{equation}
The final network output is also linear combinations of the last hidden layer
\begin{equation*}
	y = \beta\cdot y_L.
\end{equation*}
Here, $\odot$ denotes an element-wise multiplication, and
$$\theta=\{W^{(u)},b^{(u)},W^{(v)},b^{(v)},(W_l,b_l)_{l=1}^{L},\beta\}$$ are all trainable parameters.

The PirateNet architecture employed in our numerical examples starts with the same encoders and one layer of the modified MLP as above, followed by $L$ residual blocks. For $l=1, \dots, L$, the residual block $y_{l+1}=H(y_l;\theta_l,U,V)$ is defined as:
\begin{equation}
	\begin{aligned}
		f_l &= \sigma(W_l^{(1)} y_l + b_l^{(1)}),\\
		z_l^{(1)} &= f_l \odot U + (1 - f_l) \odot V,\\
		g_l &= \sigma(W_l^{(2)} z_l^{(1)} + b_l^{(2)}),\\
		z_l^{(2)} &= g_l \odot U + (1 - g_l) \odot V,\\
		h_l &= \sigma(W_l^{(3)} z_l^{(2)} + b_l^{(3)}),\\
		y_{l+1} &= \alpha_l h_l + (1 - \alpha_l) y_l,
	\end{aligned}
	\label{eq::residualblock}
\end{equation}
and the full forward pass is:
\begin{equation}
	\left\{
	\begin{aligned}
		U &= \sigma(W^{(u)}x+b^{(u)})\\
		V &= \sigma(W^{(v)}x+b^{(v)}),\\
		z & = \sigma(W_1x + b_1),\\
		y_1 & = z\odot U + (1-z)\odot V,\\
		y_{l+1} &= H(y_{l};\theta_l,U,V),\qquad  l = 1,\cdots,L,\\
		y &= \beta\cdot y_L,
	\end{aligned}\right.
	\label{st::pnet}
\end{equation}

The trainable parameters include all encoder and block components:
$$
\theta = \left\{W^{(u)}, b^{(u)}, W^{(v)}, b^{(v)}, (W_1, b_1), (W_l^{(1,2,3)}, b_l^{(1,2,3)}, \alpha_l)_{l=1}^{L}, \beta \right\}.
$$

Each residual block in \eqref{eq::residualblock} comprises three fully connected layers and two gating operations, followed by an adaptive residual connection. In practice, the coefficients $\alpha_l$ are initialized to zero, such that the network reduces to a shallow modified MLP at initialization. This design mitigates the challenges associated with training deep residual networks at initialization. The network's expressive capacity is then gradually recovered as training progresses.

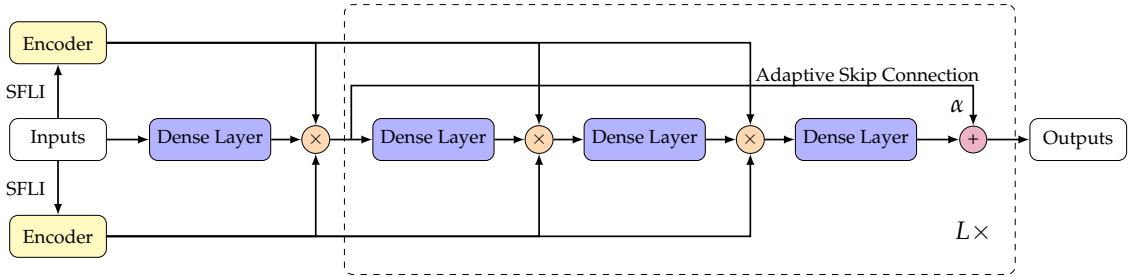
\begin{figure}[htbp]
	\centering
	\resizebox{\textwidth}{!}{  
		\begin{tikzpicture}[
			font=\footnotesize,
			node distance=0.9cm and 0.7cm,
			box/.style={rectangle, draw=black, rounded corners, minimum height=0.7cm, minimum width=1.6cm, text centered, fill=blue!30},
			box0/.style={rectangle, draw=black, rounded corners, minimum height=0.7cm, minimum width=1.6cm, text centered, fill=yellow!30},
			input/.style={rectangle, draw=black, rounded corners, minimum height=0.7cm, minimum width=1.6cm, text centered},
			mult/.style={circle, draw=black, fill=orange!30, inner sep=1.2pt},
			plus/.style={circle, draw=black, fill=purple!30, inner sep=2pt},
			skip/.style={draw=blue!50!black, thick, dashed, rounded corners},
			>=latex
			]
			
			\node[input] (input) {Inputs};
			
			\node[box0, above=of input] (encoder0) {Encoder};
			\node[box0, below=of input] (encoder1) {Encoder};
			\draw[->, thick] (input.north) -- node[left, xshift=-2pt] {SFLI} (encoder0.south);
			\draw[->, thick] (input.south) -- node[left, xshift=-2pt] {SFLI} (encoder1.north);
			
			\node[box, right=of input] (dense0) {Dense Layer};
			\draw[->, thick] (input.east) -- (dense0.west);
			
			\node[mult, right=0.5cm of dense0] (times0) {$\times$};
			\draw[->, thick] (dense0.east) -- (times0.west);
			\draw[->, thick] (encoder0.east) -| (times0.north);
			\draw[->, thick] (encoder1.east) -| (times0.south);
			
			\node[box, right=of times0] (dense1) {Dense Layer};
			\node[mult, right=0.5cm of dense1] (times1) {$\times$};
			\draw[->, thick] (times0.east) -- (dense1.west);
			\draw[->, thick] (dense1.east) -- (times1.west);
			
			\node[box, right=0.5cm of times1] (dense2) {Dense Layer};
			\node[mult, right=0.5cm of dense2] (times2) {$\times$};
			\draw[->, thick] (times1.east) -- (dense2.west);
			\draw[->, thick] (dense2.east) -- (times2.west);
			
			\node[box, right=0.5cm of times2] (dense3) {Dense Layer};
			\draw[->, thick] (times2.east) -- (dense3.west);
			
			\node[plus, right=of dense3] (plus) {+};
			\draw[->, thick] (dense3.east) -- (plus.west);
			
			\draw[->, thick] (encoder0.east) -| (times1.north);
			\draw[->, thick] (encoder1.east) -| (times1.south);
			\draw[->, thick] (encoder0.east) -| (times2.north);
			\draw[->, thick] (encoder1.east) -| (times2.south);
			\path (times0) -- (dense1) coordinate[pos=0.5] (midpoint);
			\draw[->, thick] 
			(midpoint) -- ++(0,0.9) 
			-- node[right, xshift=1.4cm, yshift=4pt] {Adaptive Skip Connection}($(plus.center)+(0,0.9)$) 
			-- node[left] {\large{$\alpha$}} (plus.north);
			
			\node[draw,dashed,rounded corners,inner sep=8pt,minimum width=11.1cm,minimum height=4.5cm,fit=(dense1) (times1) (dense2) (times2) (dense3) (plus)] (resblock) {};
			\node[below=1cm of plus] {\large{$L\times$}};
			\node[input, right=of plus] (output) {Outputs};
			\draw[->, thick] (plus.east) -- (output.west);
			
		\end{tikzpicture}
	}
	\caption{Architecture of PirateNet with Structured First-Layer Initialization (SFLI).}
	\label{fig::piratenet-architecture}
\end{figure}

To promote a high initial $\varepsilon$-rank, we apply SFLI to the encoders $U$ and $V$ defined in \eqref{eq::encoder}. Numerical evidence shows that this improves both convergence rate and final accuracy. Figure~\ref{fig::piratenet-architecture} presents a schematic diagram of PirateNet architecture with SFLI applied to its encoder layer.

\section{Numerical Experiments}\label{sec::Experiments}
This section presents a series of numerical experiments designed to evaluate the effectiveness of the proposed structured first-layer initialization (SFLI) strategy in two contexts: function approximation and the solution of PDEs using physics-informed neural networks (PINNs). Except for the first layer, all subsequent layers adopt the hyperbolic tangent activation function.
We focus on evaluating training efficiency, approximation accuracy, and the evolution of the $\varepsilon$-rank. The detailed experiment configurations and hyper-parameter settings of each example are provided in Appendix \ref{sec::appendix}.

\subsection{Function Approximations}
\begin{example}[High- and Low-Frequency 2D Function]\label{ex::fun}
	Consider the following composite target function with distinct high- and low-frequency components:

		$$f(\x) = \cos(x_1)\cos(x_2)+\cos(10x_1)\cos(10x_2), \quad \x \in [-1, 1]^2.$$
		This benchmark function is designed to test the effectiveness of the SFLI pre-training strategy under varying frequency conditions. This example is used to assess the impact of SFLI with different activation functions.  The first layer uses various SFLI pre-activations (Gaussian, Tanh, Cosine, and Hat) to evaluate their effects on $\varepsilon$-rank and convergence.
\end{example}

\begin{figure}[htbp]
	\centering
	\includegraphics[width=\textwidth]{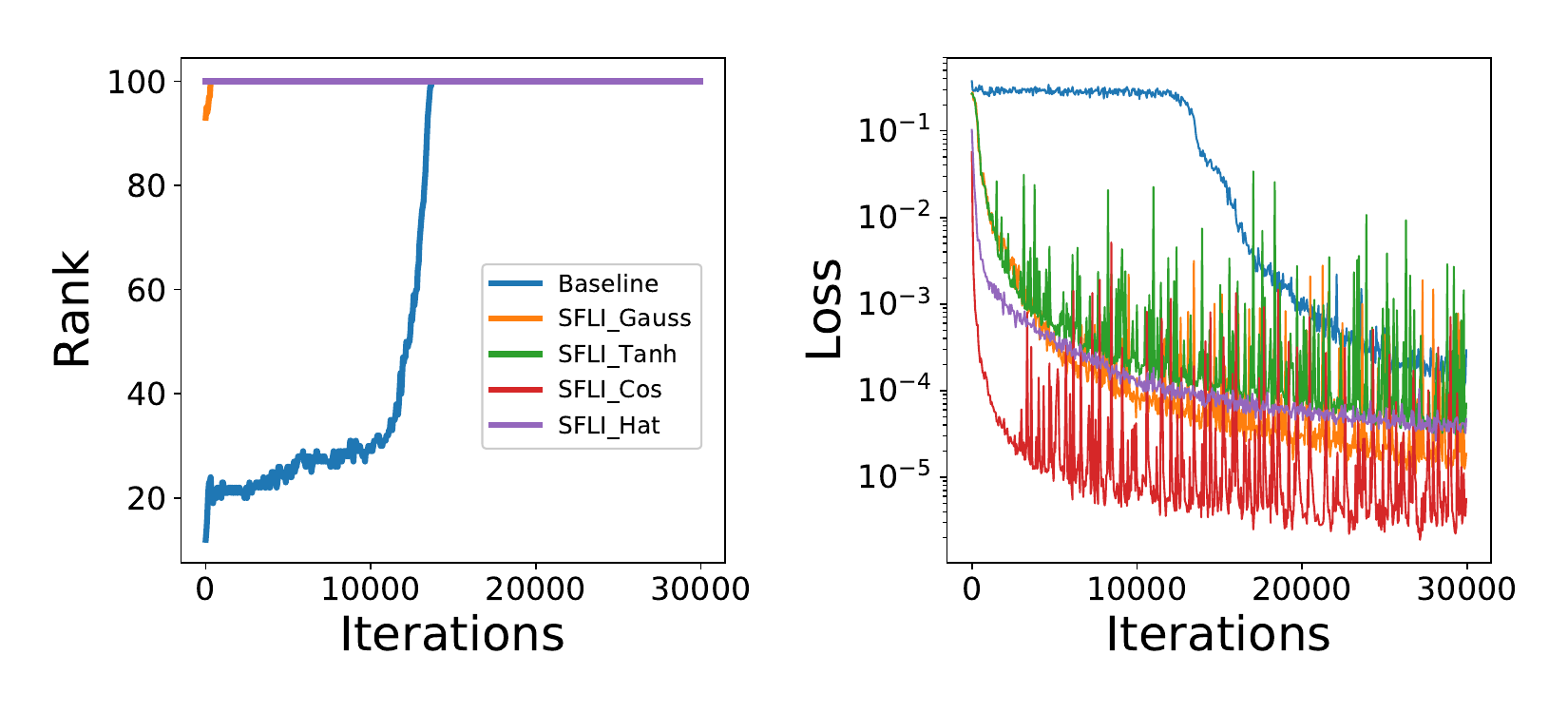}
	\caption{Comparison of $\varepsilon$-rank dynamics (left) and training loss (right) across different first-layer initializations: baseline (no SFLI) vs. SFLI with Gaussian, Tanh, Cosine, and Hat pre-activations. The detailed hyper-parameter settings are presented in Table \ref{tb::config}.}
	\label{fig::fun}
\end{figure}

As shown in Figure \ref{fig::fun}, all SFLI-based methods demonstrate superior performance compared to the baseline. Specifically, the SFLI variants exhibit high initial $\varepsilon$-ranks followed by rapid loss reduction during training, while the baseline remains trapped in a rank plateau with sluggish optimization progress. These empirical results demonstrate that the SFLI pre-training strategy effectively accelerates convergence and enhancing accuracy through $\varepsilon$-rank. Notably, the SFLI-cos method achieves the best performance in this specific case, likely due to alignment between the cosine activation and the spectral structure of the target function.

\begin{figure}[htbp]
	\centering
	\includegraphics[width=\textwidth]{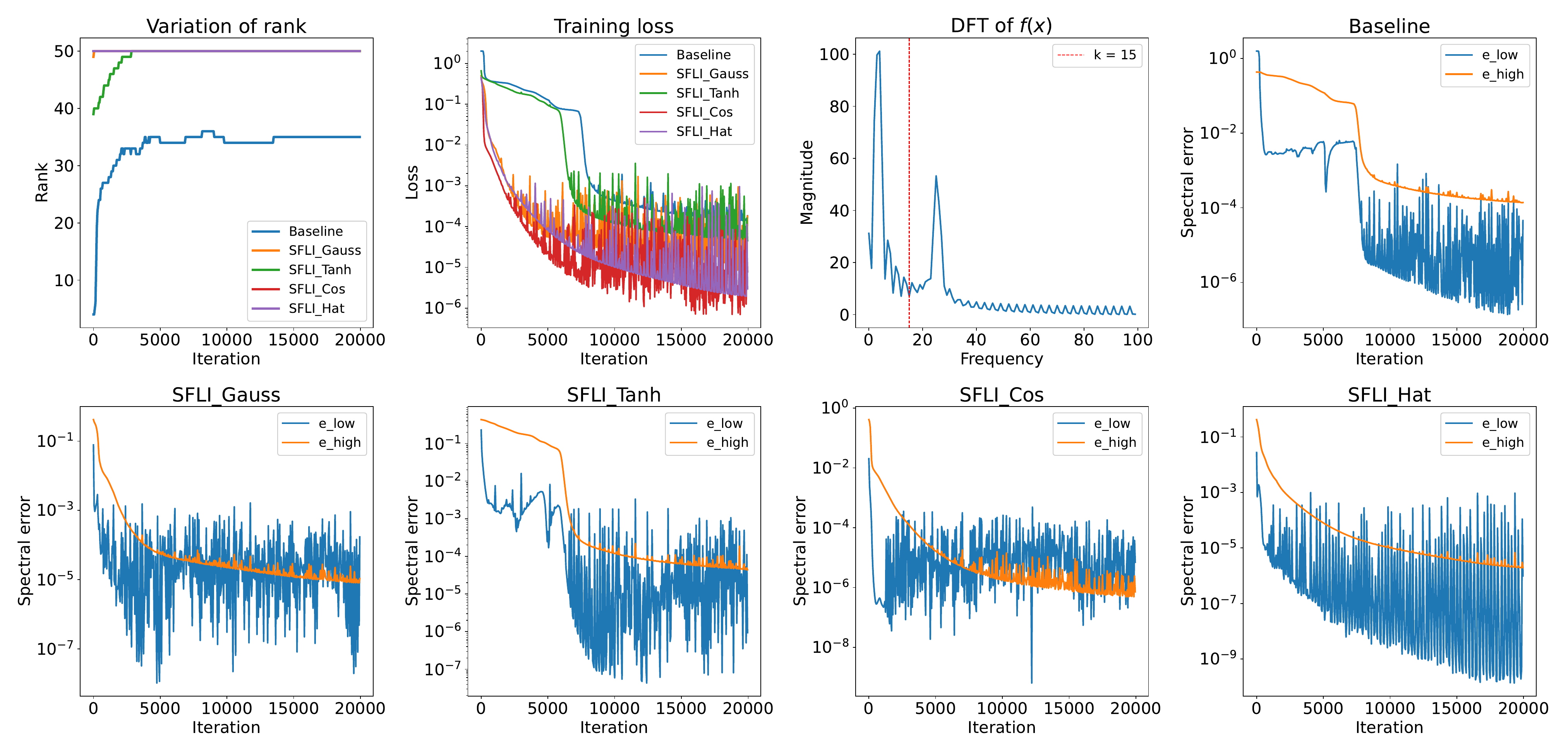}
	\caption{Spectral analysis on Example \ref{ex::spectral_bias} using different SFLI strategies. The target function contains both low- and high-frequency components. The top row displays (from left to right): (1) the evolution of $\varepsilon$-rank during training; (2) training loss curves; (3) spectrum of the target function, dividing low-frequency and high-frequency components by $k=15$; (4) the spectral fitting error of the baseline model (without SFLI), separately plotted for low- and high-frequency modes. The bottom row illustrates the evolution of spectral fitting error under different SFLI variants (Gaussian, Tanh, Cosine, Hat), where e\_low and e\_high denote the errors on low- and high-frequency components, respectively. The detailed hyper-parameter settings are presented in Table \ref{tb::config}.}
	\label{fig::frequency}
\end{figure}

\begin{example}[Discontinuous and Multiscale 1D Function]\label{ex::spectral_bias}
Consider the following piecewise-defined target function:
		$$
			f(x) =
			\begin{cases}
				(x^2 + 1)\sin(80x), & \text{if } -1 \leq x < -\frac{1}{3}, \\
				(-2x + 3)\cos(10x), & \text{if } -\frac{1}{3} \leq x < \frac{1}{3}, \\
				x^3 - x, & \text{if } \frac{1}{3} \leq x \leq 1.
			\end{cases}
		$$
		This function exhibits both discontinuities and multi-scale characteristics, making it a suitable test case for evaluating the spectral bias in neural network training. The Fourier series representation of $f$ is given by:
		\begin{equation*}
			f(x) := \sum_{k=-\infty}^\infty \hat{f}_k e^{ik\pi x}. 
		\end{equation*}
		To quantify the spectral bias, we compute the spectral error of a neural network prediction $y(x)$ in the frequency domain. The spectral error is defined as:
		\begin{equation}
			e_{\text{low}} = \sum_{k=-\delta}^\delta (\hat{y}_k - \hat{f}_k)^2, \quad e_{\text{high}} = \sum_{|k|>\delta} (\hat{y}_k - \hat{f}_k)^2,
		\end{equation}
		where $\hat{y}_k$ and $\hat{f}_k$ denote the Fourier coefficients of the predicted and target functions, respectively.
		For this particular function, the dominant frequencies are located approximately at $k=4$ and $k=25$. Therefore, we choose the spectral cutoff threshold  $\delta =15$ to separate low- and high-frequency components.

\end{example}
Figure \ref{fig::frequency} illustrates how structured first-layer initialization influences the spectral bias of neural networks. In the baseline training without SFLI, the model exhibits a pronounced spectral bias: low-frequency components are fitted quickly, while high-frequency components converge slowly and remain underfit for a long duration. In contrast, all SFLI variants significantly reduce both low- and high-frequency errors early in training. Notably, the Cosine and Gaussian based SFLI demonstrates a fast convergence in both frequency regimes, aligning with the spectral structure of the target function. These results indicate that SFLI mitigates spectral bias by encouraging $\varepsilon$-linear independence in neuron activations, allowing the network to capture fine-scale features more efficiently from the outset.


To further evaluate the effectiveness of SFLI in high dimensions, we consider a smooth function defined on the $d$-dimensional hypercube:

\begin{example}[High-Dimensional Smooth Function]\label{ex::fun_nd}
	We consider the high-dimensional function $f(\mathbf{x}) = \cos(|\mathbf{x}|^2)$ for $\mathbf{x} \in [-1,1]^d$. The performance of SFLI is tested for different dimensions $d$.
\end{example}

We compare the test relative errors after 20000 training steps obtained with and without SFLI-Gauss method across different input dimensions $d$. For SFLI, we report results using a fixed scaling factor $C=1$ as well as the best performance achieved through a grid search of $C$ over the recommended interval $[0.5, 2]$. All experiments use the same network architecture and training configurations to ensure a fair comparison. The results are summarized in Table \ref{tb::dim_comparison}, where each result is averaged over five independent runs. The detailed hyper-parameter settings are presented in Table \ref{tb::config}.

\begin{table}[htbp]
	\centering
	\caption{Test relative errors in high-dimensional function approximation.}
	\label{tb::dim_comparison}
	\renewcommand{\arraystretch}{1.2}
	\begin{tabular}{c|c|c|cc}
		\hline
		\textbf{$d$} & \textbf{Baseline} & \textbf{SFLI: $C=1$} & \multicolumn{2}{c}{\textbf{SFLI: $C=C^*$}} \\
		& Test Rel. Error   & Test Rel. Error     & $C^*$ & Test Rel. Error \\
		\hline
		5  &    \num{1.79e-2}    &    \num{2.96e-3}     &    0.9     &     \num{2.14e-3}    \\
		10 &    \num{3.50e-2}    &    \num{3.10e-3}     &    0.9     &     \num{2.86e-3}    \\
		20 &    \num{8.96e-2}    &    \num{6.89e-3}     &    1.2     &     \num{5.01e-3}    \\
		50 &    \num{9.96e-1}    &    \num{3.11e-2}     &    0.6     &     \num{1.46e-2}    \\
		\hline
	\end{tabular}
\end{table}

As shown in Table \ref{tb::dim_comparison}, the baseline model suffers from a rapid deterioration in accuracy as the input dimension increases. In contrast, SFLI significantly reduces the test error in all cases, even when using a fixed scaling factor $C=1$. Moreover, a simple grid search over $C \in [0.5, 2]$ often yields further improvement, indicating that SFLI delivers both robust default performance and enhanced accuracy with minimal tuning.

These results also validate the recommended choice of the shape parameter
\[
\gamma = C\frac{n^{\frac{1}{d}}-1}{|\Omega|^{\frac{1}{d}}},
\]
which explicitly adapts to the problem dimension $d$ and the layer width $n$, allowing SFLI to maintain stable and effective behavior across varying dimensions. The interval $C \in [0.5, 2]$ provides a practical default range for high-dimensional function approximation tasks.

\subsection{Solving Partial Differential Equations}\label{sec::pinn}
In this subsection, we demonstrate the effectiveness of the proposed structured first-layer 
initialization pre-training method in solving several benchmark partial differential equations (PDEs).
 Our experiments are based on the Physics-Informed Neural Networks (PINNs) framework~\cite{raissiPhysicsinformedNeuralNetworks2019, sirignanoDGMDeepLearning2018, moConvergenceAnalysisPinns2025}, which has achieved remarkable success in scientific computing and PDE modeling~\cite{jagtapExtendedPhysicsinformedNeural2020, JAGTAP2020113028, yingAccurateAdaptiveDeep2025, chenPhysicallyGuidedNeural2024, eskinAreTwoHidden2025, grekasDeepRitzFinite2025, badiaAdaptiveFiniteElement2025}.
Numerical examples follow the training pipeline and hyper-parameter recommendations from~\cite{wangExpertsGuideTraining2023}.
To ensure robustness and competitiveness, we incorporate several state-of-the-art PINN training techniques, including random weight factorization (RWF), periodic boundary condition embedding~\cite{dongMethodRepresentingPeriodic2021}, loss balancing~\cite{wangUnderstandingMitigatingGradient2021, wangWhenWhyPINNs2022}, and causal training~\cite{wangRespectingCausalityTraining2024}. 
We consider both the modified MLP and PirateNet architectures to demonstrate that SFLI can be effectively integrated with diverse network designs. For simplicity, we present comparisons between baseline models and those using the SFLI-Gaussian variant, which serves as a representative instantiation of our method.

Model training is conducted using mini-batch gradient descent with the Adam optimizer, where collocation points are randomly sampled at each iteration. Exact periodic boundary conditions are enforced in all numerical experiments, whenever applicable, to avoid extra loss constraints. For all examples, we compare the results with and without the application of SFLI-Gaussian, while keeping all other settings—including network architecture and hyper-parameters—identical. Detailed architecture and hyper-parameter settings of the training pipeline are provided in the Appendix \ref{sec::appendix}.

\begin{example}[Allen--Cahn equation]\label{ex::AC}
	The Allen--Cahn equation models phase separation in multi-component alloys. We consider the one-dimensional form:
	\begin{equation}
		\begin{aligned}
			&u_t - 0.0001u_{xx} + 5u^3 - 5u =0,\quad  &t \in [0,1], x \in [-1,1], \\
			&u(x, 0) = x^2\cos(\pi x), \quad  &x \in [-1,1]. \\
		\end{aligned}
	\end{equation}
\end{example}

\begin{figure}[htbp]
	\centering
	\includegraphics[width=\textwidth]{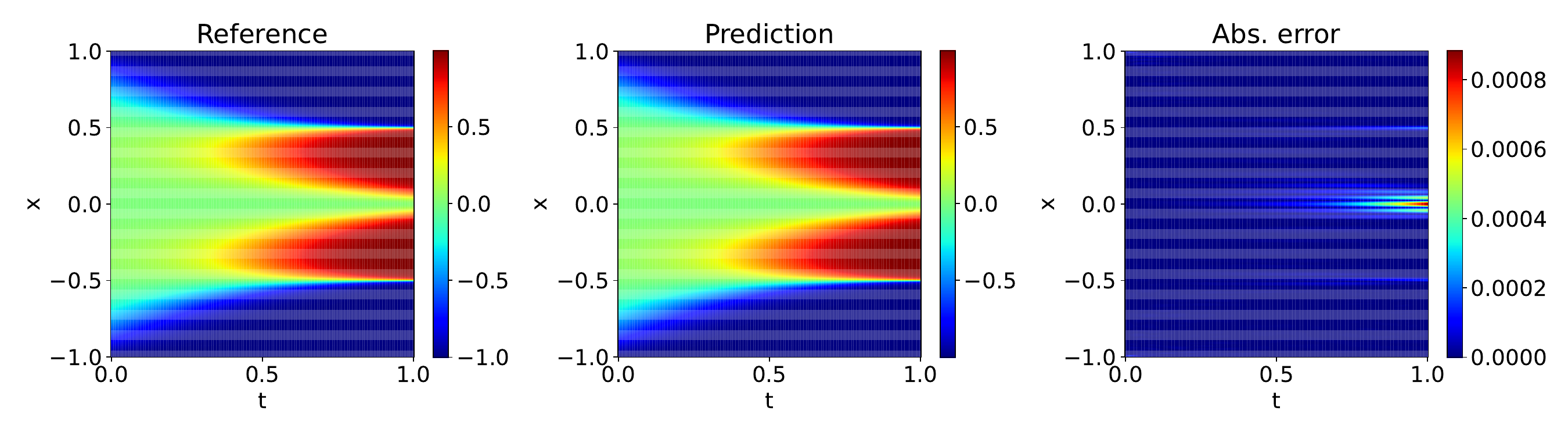}
	\caption{Allen--Cahn equation: Comparison between the solutions predicted by a trained PirateNet with SFLI and the reference solution.}
	\label{fig::AC}
\end{figure}

Figure~\ref{fig::AC} shows a comparison between the reference solution and the solution predicted by a trained PirateNet model using SFLI. The agreement between the two solutions demonstrates the model’s capability to accurately approximate the dynamics governed by the Allen--Cahn equation.

To systematically evaluate the effectiveness and generality of the proposed SFLI pre-training strategy, we consider four experimental configurations under a unified framework:
\begin{itemize}
	\item Baseline: the standard PINN model based on the PirateNet architecture, without any advanced techniques.
	\item SFLI: the Baseline model enhanced solely by the SFLI strategy.
	\item AT: the Baseline model enhanced with advanced training techniques (random weight factorization, adaptive loss weighting, and causal training).
	\item AT+SFLI: the AT model further augmented with SFLI.
\end{itemize}

All four experiments are conducted using the same PirateNet architecture and identical hyper-parameter settings, ensuring a fair comparison across different configurations. This setup allows us to isolate the effect of SFLI and assess its benefits in two complementary ways: (i) as a standalone enhancement over the standard PINN (Baseline v.s. SFLI), and (ii) as a plug-in improvement to an already advanced setting (AT v.s. AT+SFLI).

In both scenarios, the models with SFLI consistently exhibit higher $\varepsilon$-rank, faster convergence, and lower final errors, demonstrating the effectiveness of SFLI in improving both accuracy and training efficiency.

Overall, SFLI serves as a lightweight yet powerful module that can be seamlessly integrated into various PINN training pipelines to consistently boost performance.

\begin{figure}[htbp]
	\centering
	\includegraphics[width=\textwidth]{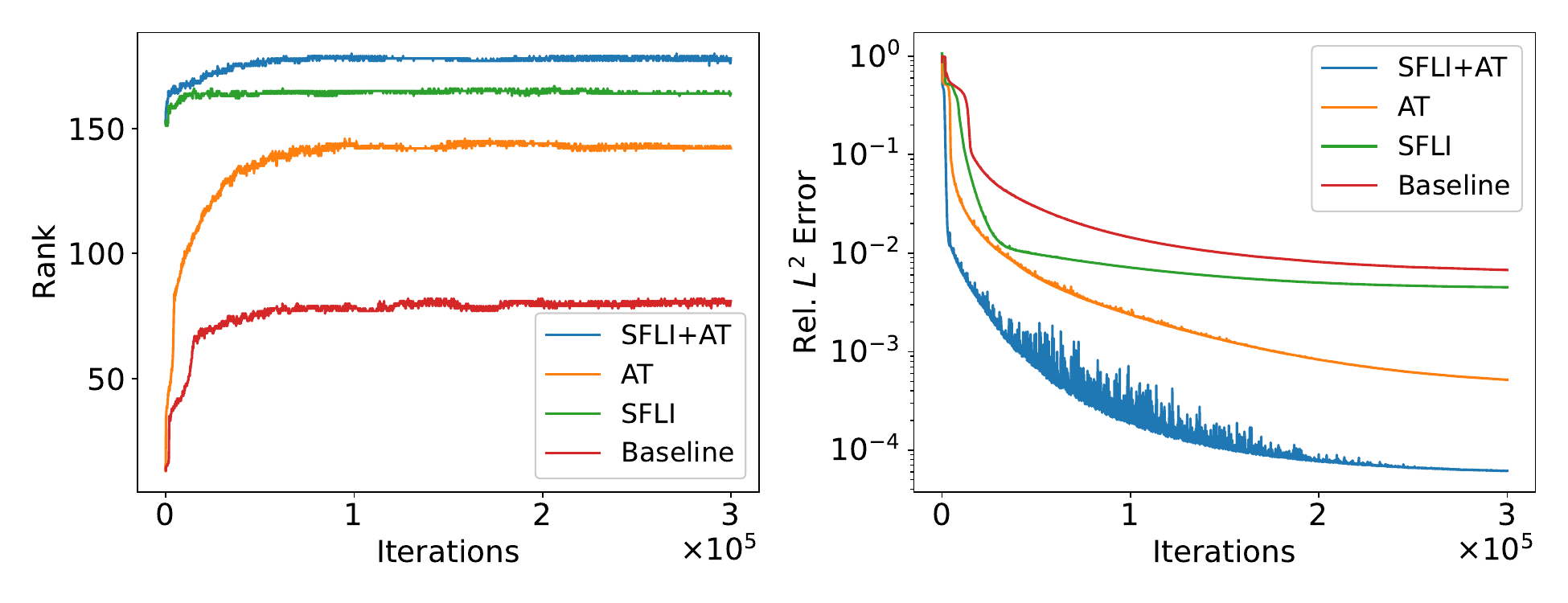}
	\caption{Allen--Cahn equation: Comparison of $\varepsilon$-rank and the relative $L^2$ errors for illustrating the impact of SFLI across different configurations, including Baseline (standard PINN), SFLI, AT (advanced techniques: RWF, loss balancing, causal training), and AT+SFLI. Detailed hyperparameter settings are provided in Table~\ref{tb::AC}.}
	\label{fig::AC_error}
\end{figure}

We now consider a classical benchmark problem in computational fluid dynamics: the lid-driven cavity flow. This problem models the steady-state motion of an incompressible viscous fluid confined within a two-dimensional square domain. The governing equations are the incompressible Navier--Stokes equations in non-dimensional form:
\begin{example}{(Lid-driven cavity flow)}\label{ex::ldc}
	\begin{equation}
		\begin{aligned}
			\bu\cdot\nabla\bu + \nabla p - \frac{1}{\text{Re}}\Delta\bu &= 0, \quad (x,y)\in (0,1)^2,\\
			\nabla \cdot \bu &= 0, \quad (x,y)\in (0,1)^2.
		\end{aligned}
	\end{equation}
	Here, $\bu=(u,v)$ denotes the velocity field in $x$ and $y$ directions, and $p$ is the pressure. The top lid moves with a constant horizontal velocity $\bu = (1,0)$, while the remaining boundaries are subject to no-slip conditions. We investigate the resulting speed field at Reynolds number $\text{Re}=3200$.
	
	To eliminate corner singularities at the top boundary, we adopt the smoothed auxiliary boundary condition from \cite{JMLR:v25:24-0313}:
	\begin{equation}\label{eq::au_boundary}
		\bu = \left(1-\frac{\cosh(50(x-0.5))}{\cosh(25)}, 0\right), \quad\text{for}\ x\in[0,1],\ y=1.
	\end{equation}
Due to the known difficulty of training PINNs at high Reynolds numbers, we employ a curriculum learning strategy to gradually increase the complexity of the problem~\cite{krishnapriyanCharacterizingPossibleFailure2021}. Specifically, the model is trained sequentially at increasing Reynolds numbers:
$\mathrm{Re} \in \{100,\ 400,\ 1000,\ 1600,\ 3200\},$
with training iterations allocated as $10^4$, $2 \times 10^4$, $5 \times 10^4$, $5 \times 10^4$, and $5 \times 10^5$ at each stage, respectively.
	
\end{example}

\begin{figure}[htbp]
	\centering
	\includegraphics[width=0.8\textwidth]{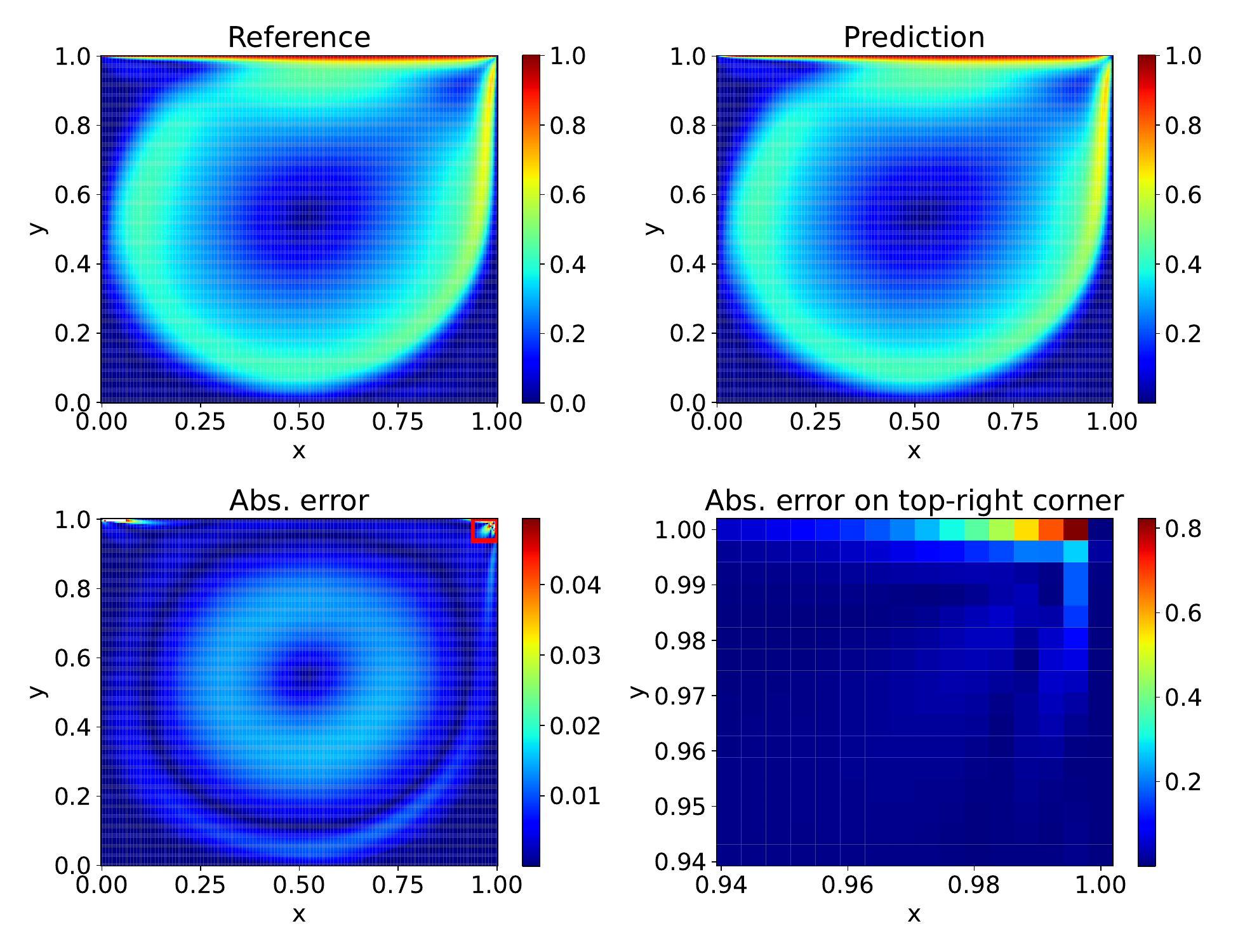}
	\caption{Lid-driven cavity flow at Re=3200: Prediction of speed $\sqrt{u^2+v^2}$ by a trained PirateNet with SFLI. The last graph is the magnified view of the top-right corner of the absolute error, corresponding to the red box in the third plot.}
	\label{fig::ldc}
\end{figure}

Figure \ref{fig::ldc} displays the predicted speed field $\sqrt{u^2 + v^2}$ for the lid-driven cavity 
flow at $\mathrm{Re} = 3200$ using PirateNet with SFLI. 
The results show close agreement with the reference solution from~\cite{ghiaHighReSolutionsIncompressible1982}, with the error primarily localized near the top corners. 

\begin{figure}[htbp]
	\centering
	\includegraphics[width=\textwidth]{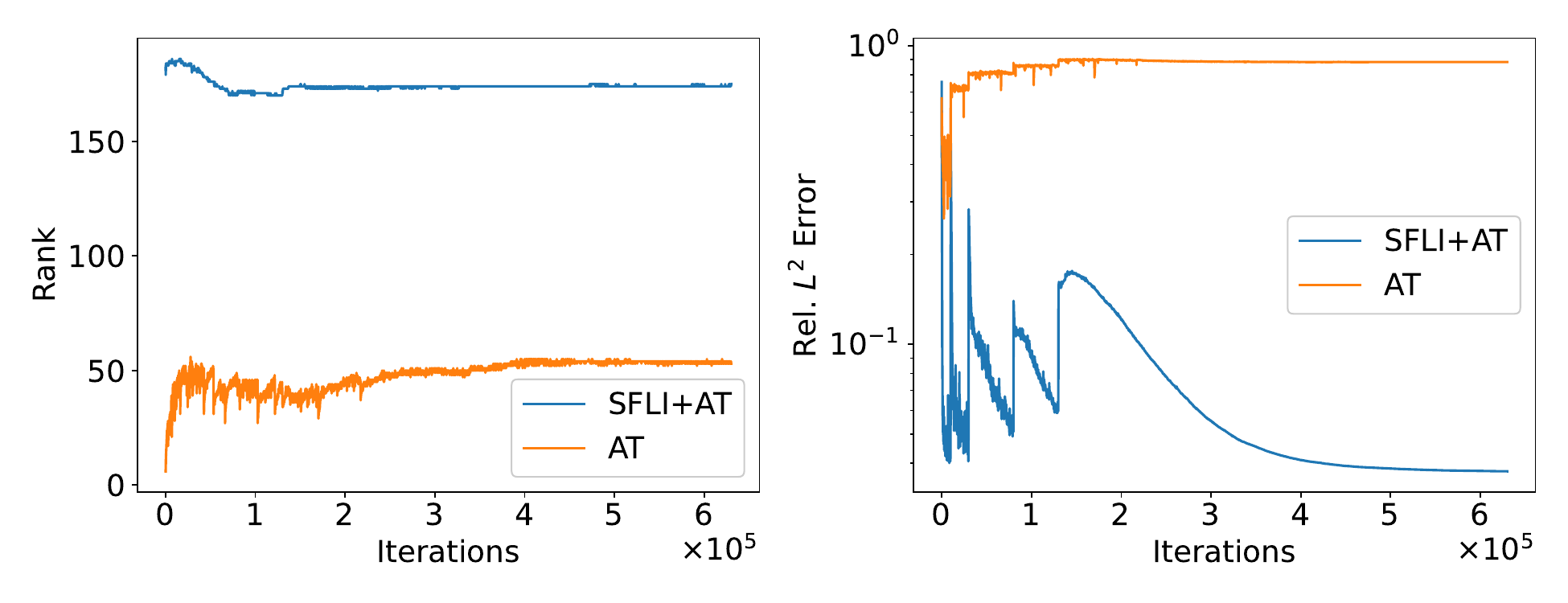}
	\caption{Lid-driven cavity flow at Re=3200: Comparison of $\varepsilon$-rank evolution and the relative $L^2$ errors of speed $\sqrt{u^2+v^2}$ between models with and without SFLI. The final relative $L^2$ errors are $\mathbf{3.75\%}$ and $88.3\%$. AT: Advanced techniques including RWF, loss balancing and curriculum training. The detailed hyper-parameter settings are presented in Table \ref{tb::ldc}.}
	\label{fig::ldc_error}
\end{figure}

To further evaluate the effectiveness of SFLI, Figure~\ref{fig::ldc_error} presents a comparison of the $\varepsilon$-rank and relative $L^2$ error between models trained with and without SFLI. 
The SFLI-enhanced model achieves a substantially higher $\varepsilon$-rank and a much lower relative error ($\mathbf{3.75\%}$ vs.\ $88.3\%$), outperforming previous baselines such as JAXPI~\cite{wangExpertsGuideTraining2023} ($15.8\%$) and the original PirateNet with random Fourier features~\cite{JMLR:v25:24-0313} ($4.21\%$). Notably, the smoothed boundary condition~\eqref{eq::au_boundary} introduces an inherent approximation error of $2.59\%$, suggesting that the observed improvement from $4.21\%$ to $3.75\%$ is a meaningful advance within the resolution limits imposed by the boundary regularization.

In the next example, we aim to demonstrate the effectiveness of our method in simulating incompressible Navier–Stokes flow  based on the velocity–vorticity formulation.
\begin{example}{(Navier–Stokes flow in a torus)}\label{ex::ns_tori}
	\begin{equation}
		\begin{aligned}
			w_t + \bu\cdot\nabla w &= \frac{1}{\text{Re}}\Delta w, \quad \text{in}\ \Omega\times[0, T] , \\
			\nabla\cdot\bu &= 0, \quad \text{in}\ \Omega\times[0, T]. \\
		\end{aligned}
	\end{equation}
	Here, $\bu = (u, v)$ represent the velocity field, $w = \nabla\times\bu$ denotes the vorticity, and Re denotes the Reynolds number. For this example, we set $\Omega=[0,2\pi]^2$, $T=10$, and $\text{Re}=100$. Periodic boundary condition is considered and the initial conditions are data from Python package \textit{jaxpi}.
	
	Our objective is to simulate the evolution of the vorticity field up to $T=10$. To this end, the temporal domain is divided into five consecutive intervals, and a time-marching strategy is employed. In each interval, a separate PINN model based on a modified MLP architecture is trained independently. The initial condition for each subsequent interval is provided by the predicted solution at the terminal time of the preceding one.
\end{example}

\begin{figure}[htbp]
	\centering
	\includegraphics[width=0.96\textwidth]{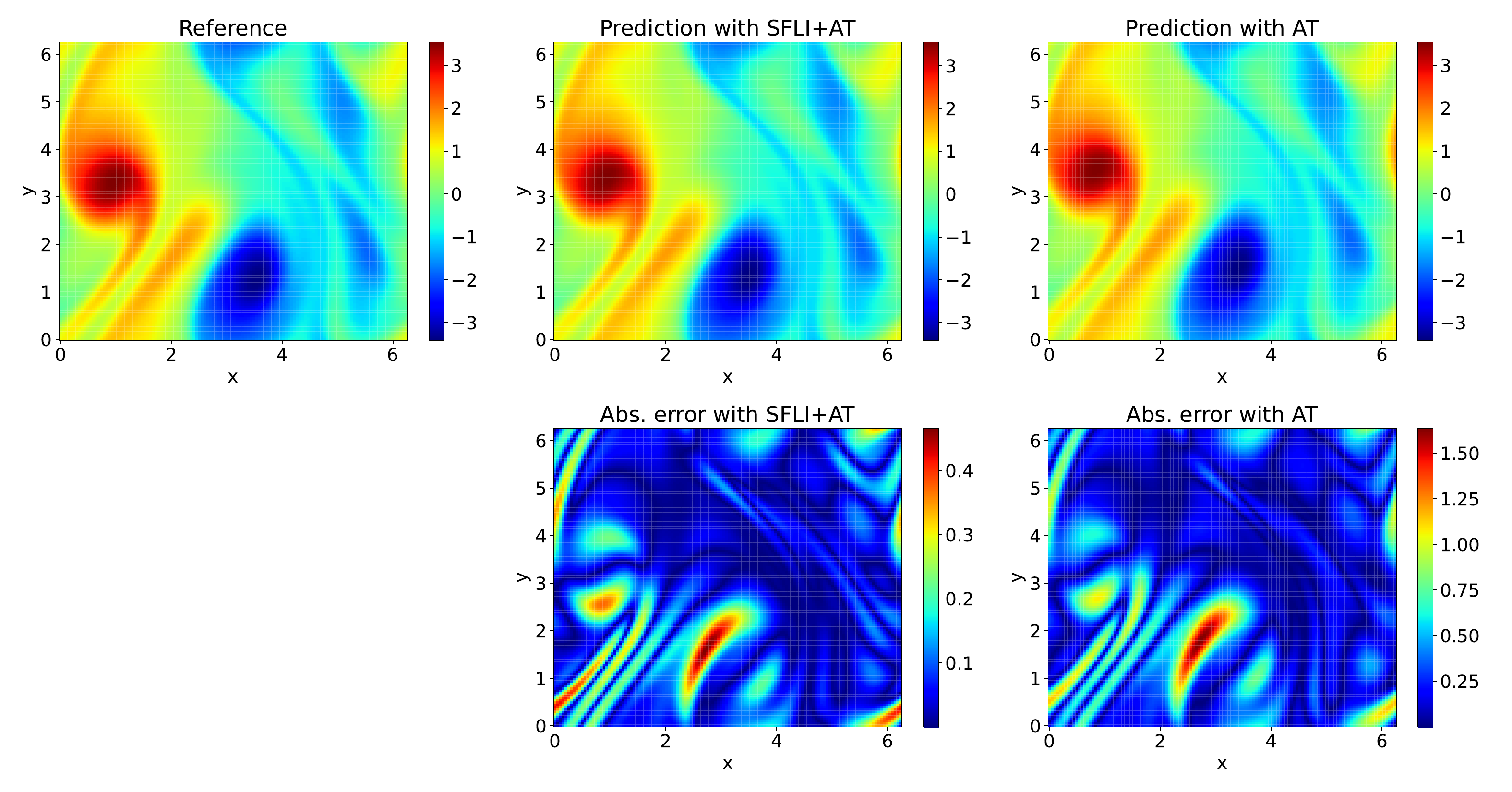}
	\caption{Navier-Stokes flow in a torus: Comparison of the predictions of $w$ at $T=10$ by a modified MLP with and without SFLI. The relative $L^2$ errors are $\mathbf{8.99\times10^{-2}}$ and $2.96\times10^{-1}$. AT: Advanced techniques including RWF, loss balancing, causal training, and time-marching strategy. The detailed hyper-parameter settings are presented in Table \ref{tb::ns_tori}.}
	\label{fig::ns_tori}
\end{figure}

Figure \ref{fig::ns_tori} compares the predicted solutions of the velocity $w$ at the final time $T=10$ for the Navier--Stokes flow in a torus, with and without the application of SFLI. The inclusion of SFLI leads to a significantly more accurate prediction, as evidenced by the lower relative $L^2$ error of \(\mathbf{8.99 \times 10^{-2}}\), compared to \(2.96 \times 10^{-1}\) in the case without SFLI. This demonstrates that SFLI enhances the neural network’s ability to capture complex vortex structures and localized flow features more accurately.

To further assess the effectiveness of SFLI in high-dimensional settings, we consider the following linear reaction-diffusion equation with a known analytical solution.

\begin{example}{(High-dimensional parabolic equation)}\label{ex::pde_nd}
	\begin{equation}
		\begin{aligned}
			u_t - \Delta u + u &= f, \quad \text{in} \ \Omega\times[0, T],\\
			u(\x, 0) &= u^*(\x, 0), \quad \text{in} \ \Omega,\\
			u(\x, t) &= u^*(\x, t), \quad \text{on} \ \partial\Omega\times[0, T].
		\end{aligned}
	\end{equation}
	Here, $\Omega=[-1,1]^d$ is the $d$-dimensional spatial domain and the final time is set to $T=0.2$. The source term is $f(\x,t) = 2d\mathrm{e}^{-t} \sin(|\x|^2) + 4 |\x|^2\mathrm{e}^{-t} \cos(|\x|^2)$. The corresponding exact solution is
	$
		u^*(\x, t) = \mathrm{e}^{-t} \cos(|\x|^2).
	$
\end{example}

We compare the relative errors of baseline PINNs and SFLI-enhanced PINNs across increasing spatial dimensions $d = 5,\ 10,\ 20,\ 50 $. For SFLI, the shape parameter is set according to the recommendation choice \eqref{choice} with a fixed scaling factor $C=1$. The results, averaged over five independent runs, are summarized in Table~\ref{tb::pde_nd}, clearly demonstrating the improved accuracy of SFLI in high-dimensional settings. The results demonstrate the improved accuracy of SFLI in high-dimensional settings. All experiments use the same network architecture and training configurations to ensure a fair comparison. Each result is averaged over five independent runs. Detailed hyper-parameter settings are provided in Table \ref{tb::parabolic}.

\begin{table}[htbp]
	\centering
	\caption{Test relative errors in solving high-dimensional PDEs.}
	\label{tb::pde_nd}
	\renewcommand{\arraystretch}{1.4}
	\begin{tabular}{c|cccc}
		\hline
		Method & \textbf{$d=5$} & \textbf{$d=10$} & \textbf{$d=20$} & \textbf{$d=50$} \\
		\hline
		\textbf{Baseline} & \num{3.12e-2} & \num{7.92e-2} & \num{1.85e-1} & \num{1.11} \\
		\textbf{SFLI}     & \num{1.73e-2} & \num{2.03e-2} & \num{3.64e-2} & \num{8.26e-2} \\
		\hline
	\end{tabular}
\end{table}

As shown in Table \ref{tb::pde_nd}, the baseline PINNs suffer a severe decline in performance as the spatial dimension increases, eventually failing to converge at $d=50$. n contrast, the SFLI-enhanced model consistently produces stable and accurate results across all tested dimensions. This highlights the robustness of SFLI in enabling effective training even in high-dimensional settings, where standard PINNs may struggle to optimize. Importantly, all results are obtained using a fixed scaling factor $C=1$, without any task-specific tuning.

\section{Discussion}

In this work, we investigated the impact of neural feature diversity on training dynamics through the lens of $\varepsilon$-rank, a quantitative measure of the effective features in neuron functions. Motivated by the insights of the staircase phenomenon, we proposed a structured first-layer initialization pre-training strategy that promotes $\varepsilon$-linear independence among neuron functions in the first hidden layer. The method is activation-function agnostic, architecture compatible, and computationally efficient, making it easy to incorporate into existing network designs. Numerical experiments on high-frequency function approximation and PDE benchmarks, including the Allen–Cahn equation, lid-driven cavity flow, Navier–Stokes systems, and high-dimensional parabolic equations, demonstrate that SFLI effectively accelerates convergence and improves predictive accuracy.

For future work, one promising direction is to incorporate $\varepsilon$-rank tracking into the training objective via dynamic regularization. Such a loss design could promote feature diversity throughout the entire training process, rather than only at initialization. Another consideration is to extend the theoretical analysis of $\varepsilon$-rank to convolutional, attention-based, or graph-based architectures. Since the current framework focuses primarily on fully connected networks, and its direct extension remains nontrivial.

In summary, this work demonstrates that structural control over the initial representational capacity of neural networks can lead to substantial improvements in training efficiency and accuracy. The proposed SFLI method offers a principled, lightweight, and broadly applicable approach for enhancing neural network performance in scientific computing tasks.

\section*{Acknowledgments}
This work is supported by the National Science Foundation of China (No.12271240, 12426312), the fund of the Guangdong Provincial Key Laboratory of Computational Science and Material Design, China (No.2019B030301001), and the Shenzhen Natural Science Fund (RCJC20210609103819018).

\bibliographystyle{plain}
\bibliography{ref}

\appendix
\section{Hyper-parameter Configurations}\label{sec::appendix}

This appendix provides the detailed hyper-parameter configurations used in all numerical experiments presented in Section \ref{sec::Experiments}. These settings include the network architecture, training strategy, optimization schedules, and shape parameters employed in the SFLI method.

Table \ref{tb::config} summarizes the experimental setups for the function fitting tasks. Tables \ref{tb::AC} to \ref{tb::parabolic} detail the complete training configurations for PDE benchmark problems, including the Allen-Cahn equation, lid-driven cavity flow, Navier-Stokes flow in a torus, and the high-dimensional parabolic equation. These configurations are chosen to ensure a fair and consistent comparison between the baseline and SFLI-enhanced models across diverse PDE types and problem scales.

\begin{table}[htbp]
	\centering
	\caption{Numerical experiment settings in function fitting examples.}
	\label{tb::config}
	\renewcommand{\arraystretch}{1.3}
	\begin{tabular}{c|c|c|c|cccc}
		\toprule
		\textbf{Example} & Layer width & \# Layers & Batch size & \multicolumn{4}{c}{\textbf{Shape parameters}} \\
		\cmidrule(lr){5-8}
		&  &  &  & Gauss & Tanh & Cos & Hat \\
		\midrule
		Example \ref{ex::fun} & 100 & 3 & 250 & 10 & 8 & 10 & 5 \\
		Example \ref{ex::spectral_bias} & 50 & 3 & 201 & 420 & 15 & 15 & 15 \\
		Example \ref{ex::fun_nd} & 128 & 3 & 1000 & Table \ref{tb::dim_comparison}\\
		\bottomrule
	\end{tabular}
\end{table}
All function approximation tasks in Table~\ref{tb::config} use multilayer perceptrons (MLPs) with \texttt{Tanh} activation unless otherwise specified. The shape parameters correspond to the localization scale $\gamma$ in the SFLI scheme (see Equation~\eqref{choice}). All models are trained using the Adam optimizer with an initial learning rate of $10^{-3}$, unless otherwise specified.

\begin{table}[htbp]
	\centering
	\caption{(Example \ref{ex::AC}) Allen-Cahn equation}
	\begin{tabular}{ll}
		\toprule
		\textbf{Parameter} & \textbf{Value} \\
		\midrule
		\addlinespace[0.3em]
		\multicolumn{2}{l}{\textbf{Architecture Parameters}} \\
		Architecture & PirateNet\\
		Number of residual blocks & 3 \\
		Layer width & 256 \\
		Activation & Tanh \\
		SFLI shape parameter & $\sqrt{10}$ \\
		Random weight factorization & $\mu = 1.0$, $\sigma = 0.1$ \\
		
		\addlinespace[0.3em]
		\multicolumn{2}{l}{\textbf{Learning rate schedule for Adam}} \\
		Initial learning rate & $10^{-3}$ \\
		Decay rate & 0.9 \\
		Decay steps & $5 \times 10^{3}$ \\
		Warmup steps & $5 \times 10^{3}$ \\
		
		\addlinespace[0.3em]
		\multicolumn{2}{l}{\textbf{Training}} \\
		Training steps & $3\times10^5$ \\
		Batch size & 8192 \\
		
		\addlinespace[0.3em]
		\multicolumn{2}{l}{\textbf{Weighting}} \\
		Weighting scheme & NTK\cite{wangWhenWhyPINNs2022}\\
		Causal tolerance & 1.0 \\
		Number of chunks & 32 \\
		\bottomrule
		\label{tb::AC}
	\end{tabular}
\end{table}

\begin{table}[htbp]
	\centering
	\caption{(Example \ref{ex::ldc}) Lid-driven cavity (Re=3200)}
	\begin{tabular}{ll}
		\toprule
		\textbf{Parameter} & \textbf{Value} \\
		\midrule
		\addlinespace[0.3em]
		\multicolumn{2}{l}{\textbf{Architecture Parameters}} \\
		Architecture & PirateNet\\
		Number of residual blocks & 6 \\
		Layer width & 256 \\
		Activation & Tanh \\
		SFLI shape parameter & 10 \\
		Random weight factorization & $\mu = 1.0$, $\sigma = 0.1$ \\
		
		\addlinespace[0.3em]
		\multicolumn{2}{l}{\textbf{Learning rate schedule for Adam}} \\
		Initial learning rate & $10^{-3}$ \\
		Decay rate & 0.9 \\
		Decay steps & $10^{4}$ \\
		Warmup steps & $5 \times 10^{3}$ \\
		
		\addlinespace[0.3em]
		\multicolumn{2}{l}{\textbf{Curriculum Training}} \\
		Reynolds Number & $[100,400,1000,1600,3200]$ \\
		Training steps & $[10^4,2\times10^4,5\times10^4,5\times10^4,5\times10^5]$ \\
		Batch size & 4096 \\
		
		\addlinespace[0.3em]
		\multicolumn{2}{l}{\textbf{Weighting}} \\
		Weighting scheme & Grad norm \cite{wangUnderstandingMitigatingGradient2021} \\
		\bottomrule
		\label{tb::ldc}
	\end{tabular}
\end{table}

\begin{table}[htbp]
	\centering
	\caption{(Example \ref{ex::ns_tori}) Navier-Stokes flow in a torus}
	\begin{tabular}{ll}
		\toprule
		\textbf{Parameter} & \textbf{Value} \\
		\midrule
		\addlinespace[0.3em]
		\multicolumn{2}{l}{\textbf{Architecture Parameters}} \\
		Architecture & Modified MLP\\
		Number of layers & 4 \\
		Layer width & 256 \\
		Activation & Tanh \\
		SFLI shape parameter & $\sqrt{1.5}$ \\
		Random weight factorization & $\mu = 1.0$, $\sigma = 0.1$ \\
		
		\addlinespace[0.3em]
		\multicolumn{2}{l}{\textbf{Learning rate schedule for Adam}} \\
		Initial learning rate & $10^{-3}$ \\
		Decay rate & 0.9 \\
		Decay steps & $2\times10^{3}$ \\
		
		\addlinespace[0.3em]
		\multicolumn{2}{l}{\textbf{Time-marching Training}} \\
		Number of time windows & 5 \\
		Training steps per window & $10^5$ \\
		Batch size & 4096 \\
		
		\addlinespace[0.3em]
		\multicolumn{2}{l}{\textbf{Weighting}} \\
		Weighting scheme & Grad norm \\
		Causal tolerance & 1.0 \\
		Number of chunks & 32 \\
		\bottomrule
		\label{tb::ns_tori}
	\end{tabular}
\end{table}

\begin{table}[htbp]
	\centering
	\caption{(Example \ref{ex::pde_nd}) High-dimensional Parabolic equations}
	\begin{tabular}{ll}
		\toprule
		\textbf{Parameter} & \textbf{Value} \\
		\midrule
		\addlinespace[0.3em]
		\multicolumn{2}{l}{\textbf{Architecture Parameters}} \\
		Architecture & MLP\\
		Number of layers & 4 \\
		Layer width & 128 \\
		Activation & Tanh \\
		SFLI shape parameter & $C=1$ in \eqref{choice} \\
		
		\addlinespace[0.3em]
		\multicolumn{2}{l}{\textbf{Learning rate schedule for Adam}} \\
		Initial learning rate & $10^{-3}$ \\
		Decay rate & 0.9 \\
		Decay steps & $2\times10^{3}$ \\
		Training steps & $2\times10^4$ \\
		
		\addlinespace[0.3em]
		\multicolumn{2}{l}{\textbf{Batch size}} \\
		Interior domain & 512 \\
		Initial domain  & 256 \\
		Boundary domain & 32$d$ \\
		
		\multicolumn{2}{l}{\textbf{Loss weight}} \\
		PDE loss & 1.0 \\
		Initial loss & 1.0 \\
		Boundary loss & 1.0 \\
		\bottomrule
		\label{tb::parabolic}
	\end{tabular}
\end{table}

\end{document}